# UNIPOTENT REPRESENTATIONS: CHANGING $q$ TO $-q$, II

P. Deligne and G. Lusztig

ABSTRACT. Consider a Chevalley group over a finite field $F_q$ such that the longest element of the Weyl group is central. In this paper we study the effect of changing $q$ to $-q$ in the polynomials which give the character values of unipotent representations of our group at semisimple elements.

## Introduction

**0.1.** Let $p$ be a prime number and let $\mathbf{k}$ be an algebraic closure of the field $F_p$ with $p$ elements. Let $G$ be connected reductive group over $\mathbf{k}$ with Weyl group $W$ and with a given maximal torus $T$. We assume that we are given an $F_p$-rational structure on $G$ with Frobenius map $F: G \to G$ such that $F(t) = t^p$ for $t \in T$ and hence $F(T) = T$.

**0.2.** We fix a prime number $l \neq p$. Let $\mathcal{R}$ be the set of all $r \in \bar{\mathbf{Q}}_l$ such that $r^2 \in \{p, p^2, p^3, \dots\}$. For $r \in \mathcal{R}$ we denote by $\mathcal{U}^r$ the set of (isomorphism classes of) unipotent representations (see [DL76, 7.8]) of the finite group $G^{F^s}$ where $r^2 = p^s$.

Let $\mathrm{Irr}(W)$ be the set of isomorphism classes of irreducible representations of $W$ over $\mathbf{Q}$. In [L84] a decomposition $\mathrm{Irr}(W) = \sqcup_c c$ into families is defined as well as a corresponding decomposition $\mathcal{U}^r = \sqcup_c \mathcal{U}^r_c$ indexed by the set $ce(W)$ of these families. In [L84] to each family $c$ is attached a finite group $\Gamma_c$, a set $M(\Gamma_c)$ (independent of $r$) and an indexing $m \mapsto \xi^r_m$ of $\mathcal{U}^r_c$ by $M(\Gamma_c)$. (we assume that a square root of $-1$ is chosen in $\bar{\mathbf{Q}}_l$).

For $\sigma \in G$ semisimple let $Z^0(\sigma)$ be the connected centralizer of $\sigma$ in $G$.

In the rest of this introduction we fix an orbit $Y$ for the conjugation action of $G$ on the set of all subgroups of $G$ of the form $Z^0(\sigma)$ for some semisimple element $\sigma \in G$. One can show that for any $s \in \{1, 2, 3, \dots\}$, the set of $G^{F^s}$-orbits (under conjugation) on the set of $F^s$-stable subgroups in $Y$ can be naturally parametrized by a finite set $\mathcal{Z}_Y$ independent of $s$ (see 1.4). In §1 to any $c \in ce(W)$, any $m \in M(\Gamma_c)$ and any $z \in \mathcal{Z}_Y$ we attach a polynomial $P_{m,z}(u) \in \mathbf{Q}[u]$ (with $u$ an indeterminate); in §2 we show that for any $s \in \{1, 2, 3, \dots\}$ and any semisimple element $\sigma \in G^{F^s}$ with $Z^0(\sigma) \in Y$ indexed by $z$ we have

(a) $\mathrm{tr}(\sigma, \xi^r_m) = P_{m,z}(p^s)$ where $r^2 = p^s$.

**0.3.** In the rest of this introduction we assume that

(a) $W$ is irreducible and the longest element $w_0$ of $W$ acts on the reflection representation of $W$ as $-1$.

Let $c \in ce(W)$ and let $m \mapsto m^!$ be the involution of $M(\Gamma_c)$ defined in [DL25]. It has the following property: if $m \in M(\Gamma_c)$ then the dimension of $\xi^r_{m^!}$ (a polynomial in $q = r^2$) is obtained (up to sign) from the dimension of $\xi^r_m$ (a polynomial in $q = r^2$) by changing $q$ to $-q$. We would like to extend this result by showing that if $\sigma$ is a semisimple element of $G^{F^s}$, $r^2 = p^s = q$, then the character of $\xi^r_{m^!}$ at $\sigma$ (a polynomial in $q$) is obtained (up to sign) from the character of $\xi^r_m$ at another semisimple element $\sigma'$ of $G^{F^s}$ (a polynomial in $q$) by changing $q$ to $-q$. But this is not true in general. (For example this fails if $G$ is of type $G_2$ and $\sigma$ has centralizer of type $A_2$.) Instead we will show that a property similar to the desired one holds for the polynomials $P_{m,z}(u)$ of 0.2. More precisely, in 1.11 we define an involution $z \mapsto z^!$ of $\mathcal{Z}_Y$ and in 1.12 we show that

(b) $P_{m^!,z^!}(u) = \pm P_{m,z}(-u)$

for any $m \in M(\Gamma_c), z \in \mathcal{Z}_Y$. In the special case where $Y$ consists of $G$ we have $z^! = z$ and (b) recovers the result of [DL25].

**0.4.** For any group $\mathcal{G}$ we denote by $cl(\mathcal{G})$ the set of conjugacy classes of $\mathcal{G}$ and by $Z_{\mathcal{G}}$ the centre of $\mathcal{G}$. For any subgroup $\mathcal{G}$ of $G$ we denote by $N\mathcal{G}$ the normalizer of $\mathcal{G}$ in $G$.

## 1. The polynomials $P_{m,z}(u)$

**1.1.** Let $Y$ be an orbit of $G$ acting by conjugation on the set of connected reductive subgroups of $G$ of the same rank as $G$. (A special case of such $Y$ was considered in 0.2.)

Let $H, H'$ be in $Y$. We can find $g \in G$ such that $H' = gHg^{-1}$. We have $NH' = gNHg^{-1}$ hence conjugation by $g$ defines an isomorphism $NH/H \to NH'/H'$ and a bijection $i_{H,H'} : cl(NH/H) \to cl(NH'/H')$ which is independent of the choice of $g$ and has an obvious transitivity property. Hence there is a well defined set $\mathcal{Z}_Y$ and bijections $j_H : cl(NH/H) \to \mathcal{Z}_Y$ for any $H \in Y$ such that $j_{H'} i_{H,H'} = j_H$ for any $H, H'$ in $Y$. This set is finite since $NH/H$ is finite for any $H \in Y$.

**1.2.** We give an alternative definition of $\mathcal{Z}_Y$. Let $Y_0 = \{H \in Y; T \subset H\}$. Note that $Y_0 \neq \emptyset$. Let $H, H'$ be in $Y_0$. We can find $g \in NT$ such that $H' = gHg^{-1}$. Then conjugation by $g$ defines an isomorphism

$$(NT \cap NH)/(NT \cap H) \to (NT \cap NH')/(NT \cap H')$$

and a bijection

(a) $\quad i'_{H,H'} : cl((NT \cap NH)/(NT \cap H)) \to cl((NT \cap NH')/(NT \cap H'))$

which is independent of the choice of $g$ and has an obvious transitivity property. Hence there is a well defined set $\mathcal{Z}'_Y$ and bijections

$$j'_H : cl((NT \cap NH)/(NT \cap H)) \to \mathcal{Z}'_Y$$

for any $H \in Y_0$ such that $j'_{H'} i'_{H,H'} = j'_H$ for any $H, H'$ in $Y_0$. It is easy to show (see 1.5(e)) that for any $H \in Y_0$, the obvious homomorphism $(NT \cap NH)/(NT \cap H) \to NH/H$ is an isomorphism. It follows that $\mathcal{Z}_Y, \mathcal{Z}'_Y$ can be canonically identified.

**1.3.** Let $H \in Y_0$. From the definition of $T$ we see that $F : G \to G$ preserves each root subgroup of $G$ with respect to $T$. Since $H$ is generated by some of these root subgroups and by $T$, we must have $F(H) = H$. It also follows that for any $H' \in Y$ we have $F(H') \in Y$.

Let $s \in \{1, 2, 3, \dots\}$. Let $Y^{F^s} = \{H' \in Y; F^s(H') = H'\}$;; this is a finite set on which the finite group $G^{F^s}$ acts by conjugation. Let $Y^{F^s}/\sim$ be the set of orbits for this action.

**Proposition 1.4.** *For any $s \in \{1, 2, 3, \dots\}$ there is a canonical bijection*

$$\mathcal{Z}_Y \to Y^{F^s}/\sim .$$

**1.5.** In this section we fix $H \in Y_0$. We define

$$\tau_H : NH \to Y^{F^s}/\sim$$

by $y \mapsto (G^{F^s} - \text{orbit of } xHx^{-1})$ where $x \in G$ satisfies $x^{-1}F^s(x) = y$ (we use Lang's theorem for $G$). This map induces a bijection

(a) $$cl'(NH) \to Y^{F^s}/\sim$$

where $cl'(NH)$ is the set of orbits of the $NH$-action on $NH$ given by $y_1 : y \mapsto y_1^{-1} y F^s(y_1)$. We show:

(b) If $y \in NH, h \in H$ then $yh = y_1^{-1} y F^s(y_1)$ for some $y_1 \in NH$.

It is enough to show that $h = y_1^{-1} y F^s(y_1) y^{-1}$ for some $y_1 \in H$. By Lang's theorem for $H$ it is enough to show that the map $H \to H$ given by $y_1 \mapsto y F^s(y_1) y^{-1}$ is a Frobenius map. By Lang's theorem for $G$ we have $y = g^{-1} F^s(g)$ for some $g \in G$. It is enough to show that the map $G \to G$ given by $y_1 \mapsto y F^s(y_1) y^{-1}$ is a Frobenius map, that is that the map $G \to G$ given by $y_1 \mapsto g^{-1} F^s(g y_1 g^{-1}) g$ is a Frobenius map. This is obvious.

We have an obvious surjective map

(c) $$cl'(NH) \to cl'(NH/H)$$

where $cl'(NH/H)$ is the set of orbits of the $NH/H$-action on $NH/H$ given by $\bar{y}_1 : \bar{y} \mapsto \bar{y}_1^{-1} \bar{y} F^s(\bar{y}_1)$. (The bijection $NH/H \to NH/H$ induced by $F^s : NH \to NH$ is denoted again by $F^s$.)

We show that

(d) The map (c) is a bijection.

It is enough to show that (c) is injective. We must show that if $y, y'$ in $NH$ and $h \in H$ are such that $hy, y'$ have the same image in $cl'(NH)$ then $y, y'$ have the same image in $cl'(NH)$. This follows from (b).

We show:

(e) The obvious homomorphism $(NT \cap NH)/T \to NH/H$ is surjective.

Let $y \in NH$. Then $yTy^{-1}$ is a maximal torus of $H$ hence $yTy^{-1} = hTh^{-1}$ for some $h \in H$ that is $y = hy'$ with $y' \in NT \cap NH$. This proves (e).

We show:

(f) $F = 1$ on $NH/H$.

By (e) it is enough to show that $F = 1$ on $(NT \cap NH)/T$. We have $(NT \cap NH)/T \subset NT/T$. But $F = 1$ on $NT/T$ since $T$ is split over $F_q$. This proves (f).

We now see that

(g) $cl'(NH/H)$ is the same as $cl(NH/H)$.

We show:

(h) The map $cl(a) : cl((NT \cap NH)/(NT \cap H)) \to cl(NH/H)$ induced by the obvious homomorphism $a : (NT \cap NH)/(NT \cap H) \to NH/H$ is a bijection.

From (e) we see that $a$ is surjective so that $cl(a)$ is surjective. We now show that $cl(a)$ is injective. Let $y, y'$ in $NT \cap NH$ be such that $hy' = xyx^{-1}$ with $h \in H, x \in NH$. By (e) we have $x = h'y_1$ where $h' \in H, y_1 \in NT \cap NH$. Hence $hy' = h'y_1 yy_1^{-1} h'^{-1}$, so that

$$y_1 yy_1^{-1} = h'^{-1} hy' h' = h'^{-1} h(y' h' y'^{-1})y' = h'' y'$$

where $h'' \in H \cap NT$. This proves injectivity of $cl(a)$. Hence $cl(a)$ is a bijection.

Combining (a),(c),(g),(h) we obtain a bijection $\beta_H : cl((NT \cap NH/(NT \cap H)) \to Y^{F^s}/\sim$.

**1.6.** We now consider $H, H'$ in $Y_0$. Let $g \in G$ be such that $H' = gHg^{-1}, T = gTg^{-1}$. Since $F = 1$ on $NT/T$ we can assume that $F(g) = g$. For $y \in NT \cap NH$ we show:

(a) $$\tau_{H'}(gyg^{-1}) = \tau_H(y).$$

(notation of 1.4.) An equivalent statement is that if $x^{-1}F^s(x) = y, x'^{-1}F^s(x') = gyg^{-1}$ then $xHx^{-1}$ and $x'H'x'^{-1} = x'gHg^{-1}x'^{-1}$ are in the same $G^{F^s}$ orbit. We have $gyg^{-1} = gx^{-1}F^s(x)g^{-1}$ hence we can take $x' = xg^{-1}$ so that $xHx^{-1} = x'H'x'^{-1}$. This proves (a).

Using (a) and the definitions we see that

(b) $\beta_{H'} i'_{H,H'} = \beta_H$.

It follows that $\mathcal{Z}_Y = \mathcal{Z}'_Y$ can be identified with $Y^{F^s}/\sim$ in such a way that $\beta_H = j'_H$ for any $H \in Y_0$. This proves 1.4.

**1.7.** If $T'$ is a torus over $\mathbf{k}$, we set $S^1_{T'} = \text{Hom}(T', \mathbf{k}^*) \otimes \mathbf{Q}$. Let $S^j_{T'}, j \in \mathbf{N}$, be the $j$-th symmetric power of $S^1_{T'}$; let $S^*_{T'} = \oplus_{j \in \mathbf{N}} S^j_{T'}$.

**1.8.** Let $H' \in Y$ and let $T'$ be a maximal torus of $H'$. Now $NT' \cap NH'$ acts naturally on the symmetric algebra $S^*_{T'/Z_{H'}}$ (preserving each summand $S^j_{T'/Z_{H'}}$). Let $\bar{S}^*_{T'/Z_{H'}}$ be the quotient of the (commutative) algebra $S^*_{T'/Z_{H'}}$ by the ideal generated by the $NT' \cap H'$-invariant elements in $\oplus_{j \in \mathbf{N}; j > 0} S^j_{T'/Z_{H'}}$, and let $\bar{S}^j(T'/Z_{H'})$ be the image of $S^j(T'/Z_{H'})$ under $S^*_{T'/Z_{H'}} \to \bar{S}^*_{T'/Z_{H'}}$. We have $\bar{S}^*_{T'/Z_{H'}} = \oplus_{j \in \mathbf{N}} \bar{S}^j_{T'/Z_{H'}}$. Note that the $NT' \cap NH'$ acts in an obvious way on $\bar{S}^j_{T'/Z_{H'}}$ for any $j \in \mathbf{N}$. Hence the space of $NT' \cap H'$-invariants $(\bar{S}^j_{T'/Z_{H'}})^{NT' \cap H'}$ is defined; it carries an action of $(NT' \cap NH')/(NT' \cap H')$.

**1.9.** Let $c \in ce(W)$; let $<,>: M(\Gamma_c) \times M(\Gamma_c) \to \bar{\mathbf{Q}}_l$ be the pairing defined in [L84, 4.14] and let $E \mapsto m_E$ be the imbedding of $c$ into $M(\Gamma_c)$ defined in [L84,§4]. Let $a_c \in \mathbf{N}, A_c \in \mathbf{N}$ be as in [DL25,2.2] and let $\Delta : M(\Gamma_c) \to \{1, -1\}$ be as in [L84,4.21]. For $E \in c$ let $b_E \in \mathbf{N}$ be as in [L84,4.1].

**1.10.** We choose a Borel subgroup $B$ of $G$ such that $T \subset B$. Then $W$ can be identified with $NT/T$ in the standard way. Hence any $E \in c$ can be regarded as a $NT/T$-module and hence as an $NT$-module.

Let $m \in M(\Gamma_c), z \in \mathcal{Z}_Y$. Let $H \in Y_0$. We define

(a) $\quad P_{m,z}(u) = \sum_{E \in c} \sum_{j \in \mathbf{N}} \mathrm{tr}(z_H, (\bar{S}^j_{T/Z_H} \otimes E)^{NT \cap H})) \Delta(m) <m, m_E> u^j \in \mathbf{Q}[u].$

Here $u$ is an indeterminate and $z_H \in (NT \cap NH)/(NT \cap H)$ corresponds to $z \in \mathcal{Z}_Y$ under the bijection $j'_H$ in 1.2. From the arguments in 1.6 we see that $P_{m,z}(u)$ does not depend on the choice of $H$.

**1.11.** In the remainder of this section we assume that $W$ satisfies 0.3(a). It follows that there exists $\dot{w}_0 \in NT$ such that $\dot{w}_0 t \dot{w}_0^{-1} = t^{-1} \mod Z_G$ for any $t \in T$. Hence for any root subgroup $U$ of $G$ with respect to $T$, $\dot{w}_0 U \dot{w}_0^{-1}$ is the root subgroup $U'$ of $G$ corresponding to minus the root corresponding to $U$. Now if $H \in Y_0$ then $H$ is generated by $T$ and by some of the root subgroups $U$ as above; moreover, if such a $U$ is contained in $H$ then so is $U'$. It follows that $\dot{w}_0 H \dot{w}_0^{-1} = H$ that is $\dot{w}_0 \in NH$. We see that

(a) $\dot{w}_0 \in NT \cap NH$ for any $H \in Y_0$.

The image of $\dot{w}_0$ in $(NT \cap NH)/(NT \cap H)$ is denoted by $\dot{w}_{0,H}$; it is an involution in the centre of $(NT \cap NH)/(NT \cap H)$ which is independent of the choice of $\dot{w}_0$. Hence multiplication by $\dot{w}_{0,H}$ gives a well defined involution of $cl((NT \cap NH)/(NT \cap H))$; under $j'_H$ this involution corresponds to an involution $z \mapsto z^!$ of $\mathcal{Z}_Y = \mathcal{Z}'_Y$, which, by arguments in 1.6, is independent of the choice of $H \in Y_0$.

**Theorem 1.12.** *Let $c \in ce(W), m \in M(\Gamma_c), z \in \mathcal{Z}_Y$. We have*

$$P_{m^!,z^!}(u) = (-1)^{A_c} P_{m,z}(-u).$$

Let $H \in Y_0$. We have

$$P_{m^!,z^!}(u) = \sum_{E \in c} \sum_{j \in \mathbf{N}} \operatorname{tr}((z^!)_H, (\bar{S}^j_{T/Z_H} \otimes E)^{NT \cap H})) \Delta(m^!) <m^!, m_E> u^j.$$

Using $\Delta(m^!) = (-1)^{a_c + A_c} \Delta(m)$, $<m^!, m_E> = (-1)^{b_E + a_c}$ (see [DL25]) we see that it is enough to show that for any $E \in c$ and any $j \in \mathbf{N}$ we have

$$\operatorname{tr}((z^!)_H, (\bar{S}^j_{T/Z_H} \otimes E)^{NT \cap H}) == (-1)^{b_E + j} \operatorname{tr}(z_H, (\bar{S}^j_{T/Z_H} \otimes E)^{NT \cap H}).$$

It is also enough to show that $\dot{w}_{0,H}$ acts on $(\bar{S}^j_{T/Z_H} \otimes E)^{NT \cap H}$ as multiplication by $(-1)^{b_E + j}$ or, more precisely, that $\dot{w}_0$ acts on $\bar{S}^j_{T/Z_H} \otimes E$ as multiplication by $(-1)^{b_E + j}$. This follows from the fact (a consequence of 0.3(a)) that $\dot{w}_0$ acts on $\bar{S}^j_{T/Z_H}$ as multiplication by $(-1)^j$ and on $E$ as multiplication by $(-1)^{b_E}$. This proves the Theorem.

## 2. Proof of 0.2(a)

**2.1.** Let $Y$ be as in 0.2. Let $s \in \{1, 2, 3, \dots\}$. We write $\Phi, q$ instead of $F^s, p^s$. Assume that $\sigma \in G^\Phi$ is semisimple and $H' := Z^0(\sigma) \in Y$ (hence $H' \in Y^\Phi$). Let $z \in \mathcal{Z}_Y = Y^\Phi/\sim$ be such that $H'$ is in the $G^\Phi$-orbit $z$. Let $\mathcal{T}_{H'}$ be the set of maximal tori of $H'$ that are $\Phi$-stable.

**2.2.** Let $H \in Y_0$. Let $T_{H'} \in \mathcal{T}_{H'}$ be such that $T_{H'}$ is as split as possible over $F_q$. We can find $x \in G$ such that $H' = xHx^{-1}$, $T_{H'} = xTx^{-1}$. We have $x^{-1}\Phi(x) \in NT \cap NH$, $\Phi(x)x^{-1} \in NT_{H'} \cap NH'$. For $t \in T_{H'}$ we have

(a) $\Phi(t) = \Phi(x)\Phi(x^{-1}tx)\Phi(x^{-1}) = \Phi(x)(x^{-1}tx)^q \Phi(x^{-1}) = \Phi(x)x^{-1}t^q x \Phi(x^{-1}).$

For $T' \in \mathcal{T}_{H'}$ we choose $g_{T'} \in H'$ such that $T' = g_{T'} T_{H'} g_{T'}^{-1}$. Note that

(b) $\Phi(g_{T'}x)(g_{T'}x)^{-1} \in NT' \cap NH'.$

Using (a) we see that $\Phi : T' \to T'$ is given by

$t' \mapsto \Phi(t') = \Phi(g_{T'})\Phi(g_{T'}^{-1} t' g_{T'})\Phi(g_{T'})^{-1} = \Phi(g_{T'})(\Phi(x)x^{-1})(g_{T'}^{-1} t' g_{T'})^q (x\Phi(x^{-1}))\Phi(g_{T'})^{-1}$

(c)
$= \Phi(g_{T'}x)(g_{T'}x)^{-1} t'^q (g_{T'}x) \Phi(g_{T'}x)^{-1}.$

**2.3.** Let $W_{H'} = (NT_{H'} \cap H')/T_{H'}$. Let $cl'(W_{H'})$ be the set of orbits of the $W_{H'}$-action on $W_{H'}$ given by $w_1 : w \mapsto w_1 w \Phi(w_1)^{-1}$ (note that $\Phi$ induces a bijection $W_{H'} \to W_{H'}$ denoted again by $\Phi$.) We define $\lambda : \mathcal{T}_{H'} \to cl'(W_{H'})$ by

$$T' \mapsto (W_{H'}\text{-orbit of image in } W_{H'} \text{ of } g_{T'}^{-1}\Phi(g_{T'})).$$

(We have $g_{T'}^{-1}\Phi(g_{T'}) \in NT_{H'} \cap H'$.) Note that $\lambda$ does not depend on the choice of $g_{T'}$. Now $H'^{\Phi}$ acts on $\mathcal{T}_{H'}$ by conjugation; let $orb(\mathcal{T}_{H'})$ be the set of orbits for this action. If $T', T''$ are in the same orbit then $\lambda(T') = \lambda(T'')$. Hence $\lambda$ induces a map

(a) $$\bar{\lambda} : orb(\mathcal{T}_{H'}) \to cl'(W_{H'}).$$

This is easily seen to be injective and using Lang's theorem for $H'$ we see that it is also surjective. Thus, (a) is a bijection. For $\mathcal{O} \in cl'(W_{H'})$ we set $\mathcal{T}_{H',\mathcal{O}} = \bar{\lambda}^{-1}(\mathcal{O})$. We have $\mathcal{T}_{H'} = \sqcup_{\mathcal{O}} \mathcal{T}_{H',\mathcal{O}}$.

**2.4.** Using 2.2(c) and Grothendieck's fixed point formula, we see that for $T' \in \mathcal{T}_{H',\mathcal{O}}$ we have

(a) $$|T'^{\Phi}| = \det(q - \Phi(g_{T'}x)(g_{T'}x)^{-1}, S_{T'}^1).$$

We show directly that the right hand side of (a) is independent of the choice of $g_{T'}$. Replacing $g_{T'}$ by $g_{T'}n$ with $n \in NT_{H'} \cap NH'$ and setting $n' = x^{-1}nx \in NT$, we have $\Phi(g_{T'}nx)(g_{T'}nx)^{-1} = \Phi(g_{T'}x)\Phi(n')n'^{-1}(g_{T'}x)^{-1} = \Phi(g_{T'}x)t(g_{T'}x)^{-1}$ (where $t \in T$ since $\Phi$ acts trivially on $NT/T$) so that $\Phi(g_{T'}nx)(g_{T'}nx)^{-1}$ acts on $S_{T'}^1$ in the same way as $\Phi(g_{T'}x)(g_{T'}x)^{-1}$. Our claim follows.

We show that, if $\mathcal{O} \in cl'(W_{H'})$ and $T' \in \mathcal{T}_{H',\mathcal{O}}$, we have

(b) $$|(NT' \cap H')^{\Phi}|/|T'^{\Phi}| = |W_{H'}|/|\mathcal{O}|.$$

We set $n = g_{T'}^{-1}\Phi(g_{T'}) \in NT_{H'} \cap H'$. Let $Z'(n)$ be the inverse image under $\pi : NT_{H'} \cap H' \to W_{H'}$ of the subgroup $Z'(w) := \{w_1 \in W_{H'}; w_1 w \Phi(w_1)^{-1} = w\}$ of $W_{H'}$ where $w = \pi(n)$. We have a well defined map $j : (NT' \cap H')^{\Phi} \to Z'(n)$ given by $\nu \mapsto g_{T'}^{-1}\nu g_{T'}$. Indeed, if $\nu \in (NT' \cap H')^{\Phi}$ we have

$$g_{T'}^{-1}\nu g_{T'} g_{T'}^{-1}\Phi(g_{T'})\Phi(g_{T'}^{-1}\nu g_{T'})^{-1} = g_{T'}^{-1}\Phi(g_{T'}).$$

Now $j$ induces a map

$$\bar{j} : (NT' \cap H')^{\Phi}/T_{H'}^{\Phi} \to Z'(n)/T_{H'} = Z'(w).$$

If $\nu \in (NT' \cap H')^{\Phi}$ satisfies $g_{T'}^{-1}\nu g_{T'} \in T_{H'}$ then $\nu \in T'^{\Phi}$. Thus $\bar{j}$ is injective. Assume now that $n_1 \in Z(n)$. Thus we have

$$n_1 g_{T'}^{-1}\Phi(g_{T'})\Phi(n_1)^{-1} = g_{T'}^{-1}\Phi(g_{T'}) \mod T_{H'}.$$

Let $\nu = g_{T'}n_1 g_{T'}^{-1}$. We have $\nu \in NT' \cap H'$ and

$$g_{T'}^{-1}\nu g_{T'} g_{T'}^{-1}\Phi(g_{T'})\Phi(g_{T'}^{-1}\nu g_{T'})^{-1} = g_{T'}^{-1}\Phi(g_{T'}) \mod T_{H'}$$

that is $\nu\Phi(\nu)^{-1} = t'$ where

$$t' \in \Phi(g_{T'})T_{H'}\Phi(g_{T'})^{-1} = \Phi(T') = T'.$$

We can write $t' = t_1'^{-1}\Phi(t_1')$ for some $t_1' \in T'$. Then we have $t_1'\nu\Phi(t_1'\nu)^{-1} = 1$, so that $t_1'\nu \in (NT' \cap H')^\Phi$ and

$$j(t_1'\nu) = g_{T'}^{-1}t_1'\nu g_{T'} = g_{T'}^{-1}t_1'g_{T'}n_1 \in T_{H'}n_1.$$

We see that $\bar{j}$ is surjective hence bijective. We deduce that

$$|NT' \cap H')^\Phi|/|T_{H'}^\Phi| = |Z'(w)|.$$

This proves (b), since $|Z'(w)| = |W_{H'}|/|\mathcal{O}|$.

**2.5.** Let $\mathcal{O} \in cl'(W_{H'})$, $T' \in \mathcal{T}_{H',\mathcal{O}}$. Let $\epsilon(H')$ (resp. $\epsilon(T')$) be $(-1)^{F_q\text{-rank of }H'}$ (resp. $(-1)^{F_q\text{-rank of }T'}$).

Using 2.4(a),(b), Borel's description [B53] of the cohomology of a flag manifold and Grothendieck's fixed point formula, we see that

$$\epsilon(H')\epsilon(T')|\mathcal{T}_{H',\mathcal{O}}| = \epsilon(H')\epsilon(T')|H'^\Phi|/|(NT' \cap H')^\Phi|$$
$$= \epsilon(H')\epsilon(T')(|H'^\Phi|/|T'^\Phi|)/(|T'^\Phi|/|(NT' \cap H')^\Phi|)$$
$$= \epsilon(H')\epsilon(T')(|H'^\Phi|/|T'^\Phi|)|\mathcal{O}|/|W_{H'}|$$

(a) $$= |H'^\Phi|_p \sum_{j \in \mathbf{N}} \text{tr}(\Phi(g_{T'}x)(g_{T'}x)^{-1}, \bar{S}^j_{T'/Z_{H'}})q^j|\mathcal{O}|/|W_{H'}|.$$

Using the isomorphism $T_{H'} \to T'$, $t \mapsto g_{T'}tg_{T'}^{-1}$, we obtain

$$\epsilon(H')\epsilon(T')|\mathcal{T}_{H',\mathcal{O}}|$$
(b) $$= |H'^\Phi|_p \sum_{j \in \mathbf{N}} \text{tr}(g_{T'}^{-1}\Phi(g_{T'})\Phi(x)x^{-1}), \bar{S}^j_{T_{H'}/Z_{H'}})q^j|\mathcal{O}|/|W_{H'}|.$$

**2.6.** Let $r \in \mathcal{R}$ be such that $r^2 = q$. From [DL76, 7.9], for $\xi \in \mathcal{U}^r$ we have

$$\text{tr}(\sigma, \xi^r)$$
$$= |H'^\Phi|_p^{-1} \sum_{T' \in \mathcal{T}_{H'}} \sum_{E \in \text{Irr}(W)} \epsilon(H')\epsilon(T')\text{tr}((g_{T'}x)^{-1}\Phi(g_{T'}x), E)(\xi : R_E)$$

where $(:), R_E$ are as in [DL25, 0.3,2.1]; note that $(g_{T'}x)^{-1}\Phi(g_{T'}x) \in NT$ acts on $E$. Now let $c \in ce(W)$ and let $m \in M(\Gamma_c)$. From [L84, 4.23] we have $(\xi_m : R_E) = \Delta(m) < m, m_E >$ if $E \in c$ and $(\xi_m : R_E) = 0$ if $E \in \text{Irr}(W) - c$. It follows that

$$\text{tr}(\sigma, \xi_m^r)$$
$$= |H'^\Phi|_p^{-1} \sum_{\mathcal{O}} \epsilon(H')\epsilon(g_{\mathcal{O}}T_{H'}g_{\mathcal{O}}^{-1})|\mathcal{T}_{H',\mathcal{O}}| \sum_{E \in c} \text{tr}((g_{\mathcal{O}}x)^{-1}\Phi(g_{\mathcal{O}}x), E)\Delta(m) < m, m_E >$$

where for each $\mathcal{O}$ we choose $g_\mathcal{O} \in H'$ such that $g_\mathcal{O} T_{H'} g_\mathcal{O}^{-1} \in \mathcal{T}_{H',\mathcal{O}}$.

Using 2.5(b) we deduce
$$\text{tr}(\sigma, \xi_m^r) = \sum_{j \in \mathbf{N}} \sum_{E \in c} \sum_\mathcal{O} \text{tr}(g_\mathcal{O}^{-1}\Phi(g_\mathcal{O})\Phi(x)x^{-1}, \bar{S}^j_{T_{H'}/Z_{H'}})$$
$$\text{tr}(g_\mathcal{O}^{-1}\Phi(g_\mathcal{O})\Phi(x)x^{-1}, \tilde{E})\Delta(m) <m, m_E> q^j |\mathcal{O}|/|W_{H'}|$$

where for any $E \in c$ we denote by $\tilde{E}$ an irreducible representation of $W_H$ obtained from $E$ using the isomorphism $NT/T \to W_{H'}$ induced by $T \to T_{H'}$, $t \mapsto xtx^{-1}$. For each $\mathcal{O}$ let $y_\mathcal{O}$ be the image of $g_\mathcal{O}^{-1}\Phi(g_\mathcal{O}) \in NT_{H'} \cap H'$ in $W_{H'}$; we have $y_\mathcal{O} \in \mathcal{O}$, see 2.3. Thus we have
$$\text{tr}(\sigma, \xi_m^r) = \sum_{j \in \mathbf{N}} \sum_{E \in c} \sum_\mathcal{O} \text{tr}(y_\mathcal{O}[\Phi(x)x^{-1}], \bar{S}^j_{T_{H'}/Z_{H'}})$$
(a) $$\text{tr}(y_\mathcal{O}[\Phi(x)x^{-1}], \tilde{E})\Delta(m) <m, m_E> q^j |\mathcal{O}|/|W_{H'}|$$

(where $[\Phi(x)x^{-1}]$ denotes the image of $\Phi(x)x^{-1} \in NT_{H'} \cap NH'$ in $(NT_{H'} \cap NH')/T_{H'}$).

We show that for $\mathcal{O}$ as above and for $v \in W_{H'}$ we have
(b) $$v y_\mathcal{O} \Phi(v)^{-1}[\Phi(x)x^{-1}] = v y_\mathcal{O}[\Phi(x)x^{-1}]v^{-1}.$$
Let $n \in NT_{H'} \cap H'$ be a representative of $v$. It is enough to show that
$$\Phi(n)^{-1}\Phi(x)x^{-1} = \Phi(x)x^{-1}n^{-1} \mod T_{H'}.$$
Since $x^{-1}nx \in NT$ and $\Phi$ acts trivially on $NT/T$ we have $\Phi(x^{-1}n^{-1}x) = x^{-1}n^{-1}xt$ for some $t \in T$. Hence we have
$$\Phi(n)^{-1}\Phi(x)x^{-1} = \Phi(x)\Phi(x^{-1}n^{-1}x)x^{-1} = \Phi(x)x^{-1}n^{-1}xtx^{-1}$$
and it remains to note that $xtx^{-1} \in T_{H'}$.

From (a),(b) we deduce
$$\text{tr}(\sigma, \xi_m^r) = \sum_{j \in \mathbf{N}} \sum_{E \in c} \sum_{y \in W_{H'}} \text{tr}(y[\Phi(x)x^{-1}], \bar{S}^j_{T_{H'}/Z_{H'}})$$
$$\text{tr}(y[\Phi(x)x^{-1}], \tilde{E})\Delta(m) <m, m_E> q^j/|W_{H'}| = \sum_{j \in \mathbf{N}} \sum_{E \in \text{Irr}(W)}$$
$$\text{tr}([\Phi(x)x^{-1}], (\bar{S}^j_{T_{H'}/Z_{H'}} \otimes \tilde{E})^{W_{H'}})\Delta(m) <m, m_E> q^j.$$

Note that $(NT_{H'} \cap NH')/T_{H'}$ acts naturally on the space of $W_{H'}$-invariants above).

Using the isomorphism $T \to T_{H'}$, $t \mapsto xtx^{-1}$, we obtain
(b) $$\text{tr}(\sigma, \xi_m^r) = \sum_{j \in \mathbf{N}} \sum_{E \in c} \text{tr}(w, (\bar{S}^j_{T/Z_H} \otimes E)^{W_H})\Delta(m) <m, m_E> q^j.$$

Here $w$ denotes the image of $x^{-1}F(x) \in NT \cap NH$ in $(NT \cap NH)/T$ (the last quotient group acts naturally on the space of $W_H = (NT \cap H)/T$-invariants above). In other words we have
$$\text{tr}(\sigma, \xi_m^r) = P_{m,z}(q),$$
see 1.10. This proves 0.2(a).

## References

[B53]  A.Borel, *Sur la cohomologie des espaces fibrés et des espaces homogenes de groupes de Lie compacts*, Ann. Math. **57** (1953), 115-207.

[DL76] P.Deligne and G.Lusztig, *Representations of reductive groups over finite fields*, Ann. Math. **103** (1976).

[DL25] P.Deligne and G.Lusztig, *Unipotent representations: changing q to $-q$*, arxiv:2508.13951.

[L84]  G.Lusztig, *Characters of reductive groups over a finite field*, Ann. Math. Studies 107, Princeton U.Press, 1984.

INSTITUTE FOR ADVANCED STUDY, PRINCETON, NJ 08540; DEPARTMENT OF MATHEMATICS, M.I.T., CAMBRIDGE, MA 02139